%
\input amstex 
\magnification=\magstep1 
\documentstyle{amsppt} 
\catcode`\@=11 
\def\logo@{} 
\catcode`\@=\active 

\baselineskip=18truept 

\NoBlackBoxes 
\TagsAsMath 
\TagsOnRight 
\parindent=.5truein 
\parskip=3pt 
\hsize=6.5truein 
\vsize=9.0truein 
\raggedbottom 
\loadmsam 
\loadbold 
\loadmsbm 
\define \x{\times} 
\let \< = \langle 
\let \> = \rangle 
\let \| = \Vert


\define \a {\alpha} 
\redefine \b {\beta} 
\redefine \dl {\delta} 
\redefine \D {\Cal D} 
\redefine \e {\epsilon} 
\define \F {\Cal F} 
\redefine \g {\gamma} 

\define \R {\bold R} 
\define \s {\sigma} 
\define \th {\theta} 
\redefine \O {\Omega} 
\redefine \o {\omega} 
\redefine \w {\omega} 
\redefine \t {\tau} 
\define \ub {\subset} 
\define \cl {\Cal L} 
\define \E {\Cal E} 
\define \A {\Cal A} 

\def \e {{\acuteaccent e}} 
\def \E {\graveaccent e} 
\def \k {\kern 5 pt} 
\def \n {\noindent} 

\redefine \i {\infty} 
\define \pt {\partial} 
\define \bl { L} 
\define \ube {\subseteq}

\let \< = \langle 
\let \> = \rangle 
\let \| = \Vert 


\redefine \b {\beta} 
\redefine \dl {\delta} 
\redefine \D {\Cal D} 

\redefine \g {\gamma} 

\redefine \O {\Omega} 
\redefine \o {\omega} 
\redefine \w {\omega} 
\redefine \t {\tau} 

\redefine \i {\infty} 


\def \n {\noindent} 
\def \r {$\R^d$} 

\def \kk {\kern 10} 
\def \k {\kern 5} 
\def \x {$X_0, \dots, X_n\,$} 
\def \y {$X_1, \dots, X_n\,$} 
\def \h {H\"ormander} 

\def \O {\Omega} 

\topmatter 

\title 
The Dirichlet problem for superdegenerate 
differential operators 
\endtitle 

\author 
Denis R\. Bell and 
Salah E.-A. Mohammed 
\endauthor 

\affil 
Mathematical Sciences Research Institute\\ 
1000 Centennial Drive\\ 
Berkeley\\ California 94720-5070. 
\medskip 
dbell\@unf.edu $\quad$ salah\@math.siu.edu 
\endaffil 

\abstract 
Let $L$ be an infinitely degenerate second-order linear operator defined 
on a bounded smooth Euclidean domain. Under weaker conditions than those 
of 
H\"ormander, we show that the Dirichlet problem 
associated with $L$ has a unique smooth classical solution. The proof 
uses the Malliavin calculus. 
At present, there 
appears to be no proof of this result using classical analytic 
techniques. 

\bigskip 

\centerline{\bf Le Probl\'eme De Dirichlet Pour } 

\centerline{ \bf Des Operateurs Differentiels Superd\'eg\'en\'er\'es} 

\noindent 
R\'E SUM\'E. Soit $L$ un op\e rateur lin\e aire d\e fini sur 
un domaine born\acuteaccent e r\e gulier de l'espace euclidien avec une 
d\e g\e n\e rescence infinie. Sous des conditions plus faibles que 
celles de 
H\"ormander, on montre que le probl\E me de Dirichlet {associ\e} 
\graveaccent a $L$ a une solution r\e guli\E re classique unique. La d\e 
monstration utilise le calcul de Malliavin. Il semble qu'il n'y ait 
\graveaccent a cette date aucune d\e monstration de ce r\e sultat par 
des 
techniques analytiques classiques. 
\endabstract 

\endtopmatter



\bigskip 

\subheading{Version fran\c caise abr\'eg\'ee} 

\medskip 
\n Soit $D$ un domaine {born\e} r\e guli\E re de $\R^d$ dont la
fronti\graveaccent
ere 
$\partial 
D$ est r\e guli\E re. Soient \break $X_0, \dots, X_n$ des champs de 
vecteurs 
et $c$ une fonction \graveaccent a valeurs r\e elles, tous d\e finis 
et r\e guliers dans un voisinage ouvert de $\bar D$. Notons $L$ 
l'op\e rateur diff\e rentiel de second ordre 
$$L = \sum_{i = 1}^n X_i^2 + X_0 + c.$$ 
Dans un article ant\e rieur [4] (voir aussi [3]) les auteurs ont 
{donn\e} des conditions suffisantes pour que $L$ soit hypoelliptique 
qui autorisent la violation sur une hypersurface de $D$ des conditions 
de H\"ormander sur l'alg\E bre de Lie. Le type de d\e g\e n\e 
rescence {autoris\e} dans [4] est {caract\e ris\e} comme suit. 

\definition{D\e finition} 

Pour chaque $k \ge 0$, soit $X^{(k)}$ la 
matrice dont les colonnes 
sont \break $X_0, \dots, X_n$, et tous les champs de vecteurs obtenus 
\graveaccent a partir 
de $X_0, \dots, X_n$ en formant les crochets de Lie it\e r\e s jusqu' 
\graveaccent a l'ordre 
$k$. Soit 
$$ \lambda^{(k)} \hbox{ la plus petite valeur propre de} 
\kern 2 pt X^{(k)}X^{(k)t}$$ 
o\graveaccent u $^t$ d\e signe la transposition des matrices. On dit 
qu'une 
hypersurface $S \subset \R^d$ est sous-critique (relativement 
\graveaccent a $L$) au point $x \in S$ quand il 
existe un 
voisinage ouvert $U$ de $x$, un entier $k \ge 0$, et un $p \in (-1, 0)$ 
tels que 
$$ 
\lambda^{(k)}(y) \ge \exp\{-[\rho(y, S)]^p\}, \quad \forall y \in U, 
$$ 
o\graveaccent u $\rho(y, S)$ est la distance euclidienne de $y$ 
\graveaccent 
a l'hypersurface S. 
\enddefinition 

Soit $K$ l'ensemble des points de $D$ o\graveaccent u $L$ ne v\e rifie 
pas la 
condition de H\"ormander. Dans 
[4], les auteurs ont {montr\e} que $L$ est hypoelliptique \graveaccent 
a condition que 
$K$ soit contenu 
dans une hypersurface $S$ de $\R^d$ de classe $C^2$ qui est 
sous-critique 
et non-carat\e ristique 
en tout point 
de $K$. En particulier, ce th\e or\graveaccent eme donne 
{l'hypoellipticit\e} des 
op\e rateurs 
$$ 
{\partial^2\over\partial x^2} + \exp(-|x|^p){\partial^2\over\partial 
y^2} + 
{\partial^2\over\partial z^2}, \quad p \in (-1, 0) 
$$ 
qui ont {\e t\e} etudi\e s par Kusuoka et Stroock dans 
[8]. Des op\e rateurs diff\e rentiels qui ont 
des d\e g\e n\e r\e scences d'ordre infini 
seront appel\e s {\it superd\e g\e n\e r\e s}. 
Le but de cet article est d'annoncer des r\e ultats qui donnent 
l'existence d'une solution 
r\e guli\E re au probl\E me de Dirichlet pour des op\e rateurs superd\e 
g\e n\e r\e s. 
Les r\e sultats sont les suivants. 

\proclaim {Th\e or\E me 1} 

Supposons que l'ensemble $K \cup \partial D$ 
est contenu 
dans une hypersurface $S$ de classe $C^2$ qui est sous-critique et
non-charact\e
ristique en 
chaque point de $K \cup \partial D$. Soient $f \in C^\infty(\bar D)$ et 
$g \in 
C^\infty(\partial D)$ et soit $u$ une solution faible (au sense des
distributions) du
probl\E me de 
Dirichlet sur 
$D$ 
$${{Lu = f \quad dans \quad D} \atop {\gamma_0 u = 
g \quad sur \quad \partial D}} \Biggr \} \eqno (1) 
$$ 
o\graveaccent u $\gamma_0$ repr\e sente la trace sectionnelle sur 
$\partial D$. Alors 
$u \in C^\infty (\bar D)$. 

\endproclaim 

\proclaim {Th\e or\E me 2} 

Supposons que les hypoth\E ses du Th\e or\E 
me 1 sont 
satisfaites et que de plus \hfil\break 
(a) $c \le 0$ sur $\bar D$ \hfil\break 
(b) Il existe $1 \le k \le d$ et $a > 0$ tels que 
$\sum_{i = 1}^n <X_i(x), e_k>^2 \ge a$ pour tout $x \in \bar D$, 
o\graveaccent u 
$e_k$ est le $k^{i\grave e me}$ vecteur de la base standard de $\R^d$. 
Alors le probl\E me de 
Dirichlet (1) admet une solution $C^\infty$ unique $u$ sur $\bar D$. 

\endproclaim 

Le Probl\E me de Dirichlet pour des op\e rateurs d\e g\e n\e r\e s $L$ a 
{\e t\e} 
{\e tudi\e} par 
divers auteurs (voir Bony [1], Derridj [6], Jerison [7], et Cattiaux 
[5]) sous 
l'hypoth\E se que $L$ satisfasse \graveaccent a la condition de 
H\"ormander en tout point. 
A la 
connaissance des auteurs, les r\e sultats ci-dessus sont les premiers 
\graveaccent a 
\e tablir l'existence d'une solution r\e guli\E re au probl\E me de 
Dirichlet sous des 
hypoth\E ses qui permettent des d\e g\e n\e r\e scences d'ordre infini. 
En particulier, 
il semble que ces r\e sultats ne peuvent \^etre obtenus \graveaccent a 
partir des 
techniques classiques. 

\medskip 

\demo{ Sch\e ma de d\e monstration de th\e or\E me 2} 

On consid\E re le processus de diffusion $d$-dimensionel {donn\e} par 
l'\e quation 
diff\e rentielle stochastique de Stratonovich 
$$ 
\left. \aligned 
d\xi^x (t) &= X_0 (\xi^x (t)) \, dt + 
\sum _{i = 1}^n X_i (\xi^x (t))\circ d W_i(t)\\ 
\xi^x(0)&=x \in D 
\endaligned \right \}\eqno(2) 
$$ 
{associ\e} \graveaccent a l'op\e rateur $L$. Alors sous les conditions 
de Th\e or\E me 2, une 
solution faible du probl\E me de Dirichlet (1) est donn\e e par 
$$ 
u(x) = E \biggl [ g(\xi^x( \tau (x))) \exp \biggl \{\int_0^{\tau(x)} 
c(\xi^x(s)\,ds \biggr \} \,dt\Big ] \kern 100 pt$$ 
$$\kern 130 pt -E \Big [ \int_0^{\tau(x)} f(\xi^x( t)) \exp \biggl 
\{\int_0^t 
c(\xi^x(s)\,ds \biggr \} \,dt \Big ] \eqno{(3)} 
$$ 
o\graveaccent u $\tau = \tau(x)$ est le premier temps de sortie de $D$ 
de la 
diffusion $\xi^x$. Sans perdre de {g\e n\e ralit\e} on peut supposer que 
$g \equiv 0$ 
sur $\partial D$, et pour simplifier prenons aussi $c \equiv 0$ ([5]). 
Donc il est 
suffisant de montrer que la fonction 
$$ 
u(x) = -E \Big [ \int_0^{\tau(x)} f(\xi^x( t))\,dt \Big ] \eqno (4) 
$$ 
est r\e guli\E re sur $\bar D$. Maintenant $L$ est hypoelliptique sur 
$D$ par 
[4] et $u$ est une solution faible de l'\e quation $Lu = 0$ dans $D$ 
([10]). 
Pour prouver la {r\e gularit\e} de $u$ jusque $\partial D$, il faut d\e 
montrer 
des estimations uniformes sur les d\e riv\e es $D^ku(x), k \ge 1$, quand 
$x$ tend 
vers $\partial D$ en restant dans $D$. Ceci est fait combinant les 
techniques de 
[2], [5], et [4]. Le r\e sultat suivant joue un role essentiel 

\proclaim {Lemme} 

Supposons que les hyoth\E ses du th\e or\E me 1 
sont satisfaites. Alors pour tout point $x_0 \in \partial D$, on peut
trouver 
un 
voisinage $U$ de $x_0$ dans $\bar D$ et un diff\e omorphisme 
$F: \bar U \to F(\bar U) \subset [0,\infty) \times \R^{d - 1}$ 
ayant les propri\e t\e s suivantes: 
\item{(i)} $F(\bar U \cap \partial D) = F(\bar U) \cap (\{ 0\} \times
\R^{d - 1})$ et 
$F(\bar U \cap D) = F(\bar U) \cap \R^{d + }$, o\graveaccent u $\R^{d +
}:=$
\hfil\break $\{x:=(x_1,x_2,\cdots,x_d): x_1 > 0 \}$. 
\item{(ii)} Il existe $1 \le i \le n$ tels que 
$F^*(X_i) = (1, 0, \dots, 0)^t$, o\graveaccent u 
$$F^*(X_i)(x) := DF(F^{-1}(x))X_i(F^{-1}(x)).$$ 
\item{(iii)} Les champs de vecteurs $F^*(X_j), j \ne i, 0 \le j \le n$,
sont de 
la forme 
$(0, Y)^t$ o\graveaccent u $Y$ est une fonction r\e guli\E re de 
$\R^d$ dans $\R^{d - 1}$. 
\item{(iv)} L'ensemble $F(\bar U \cap \partial D)$ est sous-critique 
relativement 
\graveaccent a op\e rateur {transform\e} 
$$F^*(L) = \sum_{i = 1}^n F^*(X_i)^2 + F^*(X_0) + c \circ F^{-1}.$$ 

\endproclaim 

Le lemme ci-dessus permet de r\e duire le probl\E me au case 
{o\graveaccent u} $D$ est le demi-espace 
$\bold R^{d + }$. 
En utilisant la repr\e sentation stochastique (4), on applique ensuite 
les 
estimations de [4] au processus de diffusion $F(\xi)$ pour obtenir 
les estimations d\e sir\e es sur les d\e riv\e es $D^ku(x), k \ge 1$, 
pour 
$x$ proche de $\partial D$. L'id\e e est la suivante. Comme dans [2] 
et [5], 
on conditionne la diffusion par rapport \graveaccent a son temps de 
sortie de $\bold R^{d + }$. Ceci est n\e cessaire parce que le temps de 
sortie est p.p. une 
fonction irr\e guli\E re de $x$. Les estimations de $D^ku(x), k \ge 1$, 
dont on a besoin 
sont ensuite obtenues par int\e gration par parties partielle sur 
l'espace 
de 
Wiener {conditionn\e}, comme dans le calcul de Malliavin [9]. Les 
estimations obtenues par les auteurs dans [4] sont suffisament 
robustes 
pour \^etre applicables dans le contexte pr\e sent, et jouent un role 
crucial dans l'analyse qui suit. \qed 
\enddemo 


Le th\e or\E me 1 se d\e montre comme ci-dessus. 

\medskip 

\remark{Remarque} 

Le referee a observ\acuteaccent e que l'on pouvait 
remplacer la condition b) dans le th\acuteaccent eor\graveaccent eme 2 
par l'hypoth\graveaccent ese suivante: 

\acuteaccent Pour tout $x \in \bar D$, il existe un chemin 
$\gamma: [0,\infty) 
\mapsto R^d$ tel que $\gamma(0) = x, \gamma$ quitte $\bar D$, et 
$\gamma' (t) - X_0(\gamma(t)) \in \hbox{Lie}(X_1, \dots, 
X_n)(\gamma(t))$ pour tout temps $t$ avant que $\gamma$ quitte $\bar 
D$. 

\endremark 

\bigskip

Suppose that $D$ is a bounded regular domain in $\R^d\,$ with a smooth
boundary 
$\pt 
D$. 
Assume that \x are vector fields and $c$ is a real-valued function, 
defined and smooth in an open neighborhood of 
$\bar D$. 
Let $L$ denote the second order 
differential operator 
$$L = \sum_{i = 1}^n X_i^2 + X_0 + c.$$ 
In an earlier article [4] (see also [3]) the authors gave a 
sufficient 
condition for the hypoellipticity 
of $L$ under hypotheses that allow H\"ormander's Lie algebra condition
to 
fail on 
a hypersurface in $D$. The 
type of degeneracy allowed in [4] is characterized in terms of the 
following 


\definition {Definition} 

For each $k \ge 0$, define 
$X^{(k)}$ to be a matrix with columns \x, and all 
vector fields obtained 
from \x by forming iterated Lie brackets up to order $k$. 
Define 
$$ 
\lambda ^{(k)} := \hbox{ smallest eigenvalue 
of } X^{(k)}X^{(k)t} 
$$ 
where $^t$ denotes matrix transpose. 
We say that a 
hypersurface $S \subset$ \r\kern 2 pt is 
{\it subcritical \/} (with respect 
to $L$) at $x \in S$ if there exists an open 
neighborhood $U$ of $x$, an integer $k \geq 0$, 
and $p \in (-1,0)$ such that 
$$ 
\lambda ^{(k)} (y) \geq \exp \{ - [\rho 
(y,S)]^p \}, \quad \forall 
y \in U, \eqno 
$$ 
where $\rho$ denotes the Euclidean distance between $y$ and the 
hypersurface 
$S$. 

\enddefinition 

\medskip 

Let $K$ denote the set of points in $D$ where $L$ fails to satisfy 
H\"ormander's condition. A ($C^1$) hypersurface $S$ is said to be 
{\it non-characteristic\/} (with respect to 
$L$) at $x \in K$ if at least one of the vector fields 
\y $\,$ is transversal to $S$ at $x$. 
In [4], the 
authors showed that $L$ is hypoelliptic provided $K$ is contained in a 
$C^2$ 
hypersurface $S$ of 
$\R^d$, such that $S$ is subcritical and non-characteristic at all 
points of $K$. In particular, 
this theorem asserts the hypoellipticity of the class of operators 
$$ 
{\partial^2\over\partial x^2} + \exp(-|x|^p) 
{\partial^2\over\partial y^2} + {\partial^2\over\partial z^2}, \kern 10 
pt p \in (-1, 0) 
$$ 
that were studied by Kusuoka and Stroock in [8]. 

Differential operators that exhibit infinite-order degeneracy 
will be termed 
{\it superdegenerate\/}. In this 
article we announce results that ensure 
the existence of a smooth solution to the Dirichlet problem for a large 
class of superdegenerate operators. Our results are as 
follows 


\proclaim{ Theorem 1} 

Suppose that the set $K \cup \partial D$ 
is contained in 
a $C^2$ hypersurface $S$ and 
that $S$ is subcritical and 
non-characteristic at all points of $K \cup \partial D$. Let $f \in
C^\infty (\bar 
D)$, 
$g \in C^\infty(\partial D)$ and $u$ be 
a weak solution (in the sense of distributions) to the following 
Dirichlet problem on $\bar D$ 
$${{Lu = f \quad in \quad D} \atop {\gamma_0 u = 
g \quad on \quad \partial D}} \Biggr \} \eqno (1) 
$$ 
where $\gamma_0$ denotes the sectional trace of order 0 
on $\partial D$. Then 
$u \in C^\infty(\bar D)$. 
\endproclaim 

\proclaim{Theorem 2} 

Suppose the hypotheses 
of Theorem 1 hold and in addition 
\hfil\break 
(a) $c \le 0$ on $\bar D$. \hfil\break 
(b) There exists $1 \le k \le d$ and $a > 0$ such 
that 
$\sum_{i = 1}^n <X_i(x), e_k>^2 \, \ge a, \,$for all $x \in \bar D$, 
where $e_k$ denotes the kth standard unit vector 
in \r.\hfil\break 
Then the Dirichlet problem (1) admits a unique 
smooth solution $u$ on $\bar D$. 
\endproclaim 

The Dirichlet problem for degenerate operators 
$L$ has been studied 
by several authors (cf. 
Bony [1], Derridj [6], Jerison 
[7], and Cattiaux [5]) under 
the assumption that $L$ satisfies \h's condition 
at all points. 
As far as the authors are aware, the above 
results are the first to establish 
the existence of a smooth solution to the Dirichlet 
problem 
under hypotheses that allow 
degeneracy of {\it infinite\/} order. In particular, it appears that 
these 
results are at present unavailable using classical techniques. 

\bigskip 

\demo{Outline of Proof of Theorem 2} 

Consider the $d$-dimensional diffusion process $\xi$ given by the 
Stratonovich stochastic 
differential equation 
$$ 
\left. \aligned 
d\xi (t) &= X_0 (\xi (t)) \, dt + 
\sum _{i = 1}^n X_i (\xi (t))\circ d W_i(t)\\ 
\xi(0)&=x \in D 
\endaligned \right \}\eqno(2) 
$$ 
associated with the operator $L$. 
Then by the conditions of Theorem 2, a weak solution of the Dirichlet 
problem (1) is given by 
$$ 
u(x) = E \biggl [ g(\xi^x( \tau (x))) \exp \biggl \{\int_0^{\tau(x)} 
c(\xi^x(s)\,ds \biggr \} \,dt\biggr ] 
-E \biggl [ \int_0^{\tau(x)} f(\xi^x( t)) \exp \biggl \{\int_0^t 
c(\xi^x(s)\,ds \biggr \} \,dt \biggr ] \tag{3} 
$$ 
where $\tau=\tau(x)$ is the first exit time of the diffusion $\xi^x$ 
from 
$D$. Without loss of generality, we may assume that $g \equiv 0$ on 
$\pt D$, and for simplicity we will also take $c \equiv 0$ ([5]). 
Therefore, it is sufficent to prove that the function 
$$ 
u(x) = -E \Big [ \int_0^{\tau(x)} f(\xi^x( t))\,dt \Big ] \tag{4} 
$$ 
is smooth on $\bar D$. 
Now $L$ is hypoelliptic on $D$ by [4] and $u$ is a weak solution of 
the equation $Lu=0$ in $D$ ([10]). Therefore $u$ is $C^\infty$ on $D$. 
In order to prove smoothness of $u$ up to $\partial D$, it is necessary 
to 
derive uniform bounds on the derivatives $D^{(k)} u(x), k \geq 1$, as 
$x$ 
approaches $\pt D$ through points of $D$. This is achieved by 
combining the techniques in 
[2], [5], and [4]. The following result plays a key role 

\proclaim {Lemma} 

Suppose the hypotheses of Theorem 1 hold. 
Then for each point $x_0 \in \partial D$, there exists a 
neighborhood $U$ of $x_0$ in $\bar D$ 
and a diffeomorphism $F: \bar U \to F(\bar U) \subset [0,\infty) \times
\R^{d - 1}$ 
with the following properties 
\item{(i)} $F(\bar U \cap \partial D) = 
F(\bar U) \cap (\{ 0\} \times \R^{d - 1})$ and

$F(\bar U \cap D) = F(\bar U) \cap \R^{d +}$, where $\R^{d +
}:=\{x:=(x_1,x_2,\cdots,x_d): x_1 > 0 \}$. 
\item{(ii)} There exists $1 \le i \le n$ such that 
$F^*(X_i) = (1, 0, \dots, 0)^t$, where 
$$F^*(X_i)(x):= DF(F^{-1}(x))X_i (F^{-1}(x)).$$ 
\item{(iii)} The vector fields $F^*(X_j), j \ne i, 0 \le j \le n$, 
have the form 
$(0, Y)^t$ where $Y$ is a 
smooth function from $\R^d$ to $\R^{d - 1}.$ 
\item{(iv)} The set $F(\bar U \cap \partial D)$ is subcritical with
respect 
to the 
transformed operator 
$$F^*(L) = \sum_{i = 1}^n F^*(X_i)^2 + F^*(X_0) + c \circ F^{-1}.$$ 
\endproclaim 

\n 

The above lemma reduces the problem to the case when $D$ is the 
half-space 
$\R^{d + }$. 
Using the stochastic representation (4), the estimates in [4] are then 
applied to the diffusion process $F(\xi)$ 
in order to obtain the desired estimates on $D^{(k)} u(x), k \geq 1$, 
at points close to $\partial D$. 
The idea is as follows. Following [2] 
and [5], we condition 
the diffusion $F(\xi)$ 
up to its exit time from $\R^{d + }$. This is necessary because the exit 
time is a.s. irregular as 
a function of $x$. The required estimates on $D^{(k)} u(x), k \geq 1$, 
are then obtained through a process of partial integration by parts on
the 
conditioned Wiener space, 
in the manner of the Malliavin calculus [9]. The estimates obtained by 
the authors in [4] are sufficiently robust 
to be applied in the present context, and 
play a crucial role in the subsequent analysis. 
\qed 

\enddemo 

Theorem 1 is proved by an argument similar to the above. 

\remark{Remark} 

The referee has observed that condition (b) in Theorem 2 may be replaced
by the
following hypothesis: 

For each $x \in \bar D$, there is a $C^1$ path $\gamma: [0,\infty) \to
\R^d$ such
that 
$\gamma(0)=x$, $\gamma$ leaves $\bar D$, and $\gamma'(t) -X_0(\gamma(t))
\in
\hbox{Lie}(X_1, \cdots, X_n)(\gamma(t))$ for all times $t$ before
$\gamma$ exits
$\bar D$. 

\endremark 

\subheading{Acknowledgement} 

The authors are grateful to Guy David for the translation into French,
and to the 
referee for helpful suggestions. 
The research of Denis Bell is supported in part by NSF grant 
DMS-9703852 and by MSRI, Berkeley, California. 
The research of Salah Mohammed is supported in part by NSF grants 
DMS-9503702, DMS-9703596, and by MSRI, Berkeley, California. 

\bigskip 

\bigskip 
\bigskip 

\centerline{\bf{ References}} 

\bigskip 
\baselineskip=14truept 
\parskip=3truept 

\medskip 
\item{[1]\quad} Bony J. M., Principe du maximum, in\'egalit\'e 
de Harnack et unicit\'e du probl\`eme de Cauchy pour les op\'erateurs 
elliptiques d\'eg\'en\'er\'es, {\it Ann. Inst. Fourier\/}, Grenoble, 
t. 19 (1969) 277-304. 
\medskip 

\item{[2]\quad} Ben-Arous G., 
Kusuoka S., Stroock D. W., 
The Poisson kernel for certain degenerate 
elliptic operators, {\it J. Funct. Anal.\/} 56, no. 2 (1984) 171--209. 
\medskip 

\item{[3]\quad} Bell D. R., Mohammed 
S.-E. A., 
Hypoelliptic parabolic operators with exponential 
degeneracies, {\it 
C. R. Acad. Sci. Paris\/}, t. 317, S\'erie I 
(1993) 1059-1064. 
\medskip 

\item{[4]\quad} Bell D. R., 
Mohammed S.-E. A., An extension of 
H\"ormander's theorem for infinitely degenerate 
second-order operators, {\it 
Duke Math J.\/}, 78, no. 3 (1995) 453-475. 

\medskip 
\item{[5]\quad} Cattiaux P., Calcul stochastique et op\'erateurs 
d\'eg\'ener\'es 
du second ordre II. Problem de Dirichlet, {\it Bull. Sc. math., 2$^e$ 
s\'erie\/} 
115 (1991) 
81-122. 

\medskip 
\item{[6]\quad} Derridj M., Un probl\`eme aux limites pour une classe 
d'op\'erateurs du second ordre hypoelliptiques {\it Ann. Inst. 
Fourier\/}, 
Grenoble, t. 21 (1971) 99-148. 
\medskip 

\medskip 

\item{[7] \quad} Jerison D., The Dirichlet problem for the Kohn 
Laplacian 
on the Heisenberg group I, {\it J. Funct. Anal.\/} 43, no. 1 (1981) 
97-142. 

\medskip 
\item{[8] \quad} Kusuoka S., 
Stroock D. W., Applications of the 
Malliavin calculus, Part II, {\it Journal 
of Faculty of Science, 
University of Tokyo}, Sec. 1A, 32 
(1985) 1-76. 

\medskip 
\item{[9] \quad} Malliavin P., Stochastic 
calculus of variations 
and hypoelliptic operators, {\it Proceedings 
of the International 
Conference on Stochastic Differential 
Equations, Kyoto}, 
Kinokuniya (1976) 195-263. 

\medskip 
\item{[10] \quad} Stroock D. W., Varadhan 
S. R. S., On degenerate 
elliptic-parabolic operators of second order 
and their associated diffusions, 
{\it Comm. Pure. Appl. Math.\/} 25 (1972) 
651--713. 

\end 
\enddocument